\documentclass[11pt,a4paper,reqno]{amsart}
\usepackage[foot]{amsaddr}
\usepackage{geometry}
\usepackage{algorithm}
\usepackage{algpseudocode}
\usepackage{float}
\usepackage{alphabeta}
\usepackage{minitoc}
\usepackage{listings}
\usepackage[dvipsnames]{xcolor}
\usepackage{enumerate}
\usepackage{bm}
\usepackage[utf8]{inputenc}
\usepackage[T1]{fontenc}
\usepackage{graphicx}
\usepackage{longtable}
\usepackage{wrapfig}
\usepackage{rotating}
\usepackage[normalem]{ulem}
\usepackage{amsmath}
\usepackage{amssymb}
\usepackage{capt-of}
\usepackage{hyperref}
\usepackage{textcase}

\usepackage{enumitem}
\usepackage{parskip}
\setlist[itemize,1]{leftmargin=\dimexpr 26pt-0.025in}
\usepackage{fancyhdr}
\usepackage{caption}
\theoremstyle{definition}

\renewcommand{\tableofcontents}{}

\usepackage[backend=biber, style=numeric, sorting=none,
maxbibnames=999, giveninits=true]{biblatex}

\author[C. Smaragdakis]{Costas Smaragdakis}

\address{Department of Statistics and Actuarial-Financial Mathematics,
  University of the Aegean, Samos, 83200 Karlovassi, Greece
\& Institute of Applied and Computational Mathematics, FORTH, 70013 Heraklion, Greece
}
\email{kesmarag@aegean.gr}

\title{Learning Geometric-Aware Quadrature Rules\\ for Functional Minimization}



\hypersetup{
 pdfauthor={Notes},
 pdftitle={ML Quadrature},
 colorlinks=true,
 linkcolor={white!30!black},
 citecolor={blue!50!black},
 urlcolor={blue!80!black} }
\begin{document}

\maketitle
\tableofcontents

\begin{abstract}
Accurate numerical integration over non-uniform point clouds is a challenge for modern mesh-free machine learning solvers for partial
differential equations (PDEs) using variational principles.
While standard Monte Carlo (MC) methods are not capable of handling a non-uniform point cloud, modern neural network architectures can deal with permutation-invariant inputs, creating quadrature rules for any point cloud.
In this work, we introduce \texttt{QuadrANN}, a Graph Neural Network (GNN) architecture designed to learn optimal quadrature weights directly from the underlying geometry of point clouds.
The design of the model exploits a deep message-passing scheme where the initial layer encodes rich local geometric features from absolute and relative positions as well as an explicit local density measure. In contrast, the following layers incorporate a global context vector. These architectural choices allow the \texttt{QuadrANN} to generate a data-driven quadrature rule that is permutation-invariant and adaptive to both
local point density and the overall domain shape. We test our
methodology on a series of challenging test cases, including integration on convex and non-convex domains and estimating the solution of the Heat and Fokker-Planck equations.
Across all the tests, \texttt{QuadrANN} reduces the variance of the
integral estimation compared to standard quasi-Monte Carlo methods by warping the point clouds to be more dense in critical areas where the integrands present certain singularities.
This enhanced stability in critical areas of the domain at hand is critical for the optimization of energy functionals, leading to improved deep learning-based variational solvers.
\end{abstract}

\section{Introduction}
Variational principles are fundamental to the mathematics and physical
sciences, creating the ground for solving challenging problems in
fields covering classical and statistical mechanics to modern
probability theory.
These principles often require the minimization of
a functional, 
which typically involves an integral over a particular domain $\Omega$:
\begin{equation}
  \label{eq:ju}
  \mathcal J[u] = \int_\Omega L(x,u(x),\nabla u(x))dx,\quad u\in H^1(\Omega),
\end{equation}
where $L$ is a given operator. The objective is to find a function
$u^*\in H^1(\Omega)$ that minimizes this functional, i.e., $\mathcal J[u^*] \leq
\mathcal J[u]$ for any function $u\in H^1(\Omega)$.
Such a problem is equivalent to solving the corresponding Euler-Lagrange partial differential equation (PDE), and the pursuit of its solution is a cornerstone of computational science.

For many years, the numerical solution of these problems has been
driven by mesh-based methods like the Finite Element Method
(FEM) \cite{zienkiewicz2005finite} { or global
  pseudospectral methods \cite{fornberg1994review}}.
{ While these methods are powerful, they heavily rely on domain meshing,} which can be computationally restrictive for problems in high-dimensional spaces or with complex geometries. 
In recent years, we have witnessed a shift in the methodologies that are used to solve problems, thanks to the introduction of deep learning-based solvers. 
Methods spanning the Deep Ritz Method (DRM) \cite{e2018deep} and Physics-Informed Neural Networks (PINNs) \cite{raissi2019physics} have been introduced as mesh-free alternatives. 
These modern approaches reformulate the variational problem as an
optimisation problem, where now the parameters of an Artificial Neural
Network (ANN), representing the solution $u$, are trained to minimize
the functional $\mathcal{J}[u]$.

A significant aspect of deep learning solvers is the accurate
numerical approximation of the integral shown in Equation
(\ref{eq:ju}).
The solving process depends heavily on the minimization of the discretized functional, $\hat{\mathcal{J}}[u]$, typically evaluated with respect to a given set of sample points. The accuracy of this numerical integration is critical for obtaining an accurate solution. As this work will underscore, unless the approximation $\hat{\mathcal{J}}[u]$ is highly accurate, the minimizer $\hat{u}^*$ obtained from training the ANN will generally not be a reasonable estimate of the actual minimizer $u^*$ of the original functional.

{
The simplest approach for approximating such integrals is the Monte
Carlo (MC) method. However, MC suffers from slow convergence, with an
error rate of $\mathcal{O}(n^{-1/2})$. While variants like quasi-Monte Carlo
(QMC) Methods offer faster convergence for certain classes of
integrands \cite{niederreiter1992random}, the specific requirements of
modern Deep Learning-based PDE solvers, namely, the need to integrate
over arbitrary, unstructured, and adaptive point clouds, exposing significant gaps in the existing literature.

First, well-established quadrature rules, such as Gaussian or Clenshaw-Curtis quadrature, and their sparse grid extensions (e.g., Smolyak rules) \cite{bungartz2004sparse, narayan2014adaptive, feinberg2015chaospy}, provide optimal convergence but fundamentally rely on grids. Even extensions designed to handle non-convex domains via convex optimization \cite{ryu2015extensions} typically require specific constraints or a priori knowledge of the basis, limiting their use in mesh-free, deep learning settings.

Second, for generic unstructured point clouds, Kernel Optimal Quadrature (KOQ) is often used as a benchmark. This approach determines the optimal weights by solving a linear system based on the pairwise kernel evaluations of the points \cite{belhadji2019kernel, oates2019probabilistic}. However, calculating these weights requires solving a dense linear system ($\mathcal{O}(n^3)$). In a deep learning training loop, where the point cloud changes at every optimization step, solving a linear system iteratively is computationally prohibitive.

Third, methods for minimizing discrepancy primarily focus on optimizing the geometry of the point set itself
 \cite{kirk2025low, chahine2025improving, clement2025optimization}, these
approaches aim to set the points at optimal locations.
In many physics-informed learning scenarios, the point distribution is
fixed by external factors, such as the geometry of the domain or
importance sampling based on the physics of the problem.

To address these limitations, we introduce \texttt{QuadrANN}, a Graph Neural Network (GNN) architecture that learns optimal quadrature weights directly from the underlying geometry of point clouds. Our work builds upon the growing interest in machine learning for numerical integration \cite{briol2019probabilistic, zepeda2021deep, rivera2022quadrature}, but specifically targets the gap of handling fixed, non-uniform distributions efficiently. \texttt{QuadrANN} learns a continuous mapping from arbitrary geometry to quadrature weights, enabling accurate integration for the dynamic, unstructured point clouds inherent to neural PDE solvers without the computational bottleneck of kernel methods.
}

The remainder of this work is as follows. In Section 2, we recap the neural quadrature problem for functional integrals.
Section 3 introduces our approach for generating non-uniform point clouds model for model training and testing.
In Section 4, we present the proposed QuadrANN architecture and analyze our architectural decisions.
Section 5 presents applications to demonstrate the performance of
QuadrANN, and finally, Section 6 provides the conclusion.

\section{Neural Quadrature for Functional Integrals}

Let domain $\Omega\subset \mathbb R^d$. The core idea in a quadrature rule is to estimate an integral using a weighted sum over a set of $n$ points $\mathcal{X} = \{x_1, x_2, \dots, x_n\}$ sampled from $\Omega$:
\begin{equation}
  \int_\Omega L(x,u(x),\nabla u(x))dx  \approx \sum_{i=1}^N w_i L(x_i,u(x_i),\nabla u(x_i)).
\end{equation}
The main challenge here is to find the optimal weights $(w_1, \dots, w_n)$ for a given set of points $\mathcal{X}$. These weights should depend on the overall geometry of the point cloud but be independent of the ordering in which the points are provided to a model. To achieve this, we represent the mapping from the point set $\mathcal{X}$ to the weights $(w_i)_{i=1}^n$ by means of an Artificial Neural Network (ANN) that is insensitive to the permutation of its inputs. We denote this network as $W(\cdot\ ;\theta)$, which takes the entire set $\mathcal{X}$ as input and produces $N$ corresponding weights. The integral approximation is then:

\begin{equation}
\int_\Omega L(x,u(x),\nabla u(x))dx  \approx \sum_{i=1}^n W_i(\mathcal X;\theta) L(x_i,u(x_i),\nabla u(x_i)),
\end{equation}
where $W_i(\cdot\ ;\theta)$ is the $i$-th output of the network corresponding to
the point $x_i$. While the underlying framework is applicable to more
general domains, the theoretical development of our ANN quadrature
rule will henceforth be mainly focused on the $d$-dimensional
hypercube $\Omega = [0,1]^d$.
However, we will expand our analysis to include more complex
geometries, such as a non-convex L-shaped domain in 2D (see Subsection \ref{ex2:lshaped}).

The network is trained to produce model parameters that guarantee a good
integration performance for a predefined basis set of test functions. This
basis is designed to be highly versatile, combining polynomial and
trigonometric functions to effectively approximate a broad class of
integrands. The core of the basis consists of multivariate
  polynomials of total degree at most $J$, which is fixed at $5$
throughout this work, formed by tensor products of normalized Hermite polynomials, $\{\Phi_{\boldsymbol{j}}(\cdot)\}$:
\begin{equation}
   \Phi_{\boldsymbol j}(x) = \prod_{k=1}^d \tilde H_{e_{j_k}}(x_k),\ \text{where} \sum_{k=1}^d j_k \leq J.
\end{equation}
The polynomials $\tilde H_{e_{j_k}}$ are normalised such that their integral over $[0,1]$ is unity. To increase the capability of the basis for oscillatory functions, the degree-zero Hermite polynomial is reformulated as a trigonometric term:
\begin{equation}
   \tilde H_{e_0}(x;\kappa,\varphi) = 1 + \cos(\kappa\pi x + \varphi) - \int_{\Omega} \cos(\kappa\pi x + \varphi)dx,
\end{equation}
where $\kappa$ is an integer randomly sampled from $\{1,2,\dots,
K_{\mathrm{max}}\}$ during training and \(\varphi\in[0,2\pi)\) is a random phase
.
This hybrid construction, controlled by hyperparameters $J$ and
$K_{\mathrm{max}}$, creates a powerful approximation space. The
polynomial components are well-suited for capturing smooth,
low-frequency behaviour, while the randomized trigonometric terms are
explicitly included to represent oscillatory patterns that are
difficult to approximate efficiently with a pure polynomial basis. It
is important to emphasize that this reformulation does not affect the
ability of the model to integrate constants. This property is instead
enforced by the architecture of the network, as the final \texttt{softmax} activation layer constrains the output weights to sum to unity ($\sum W_i = 1$).

The training objective is operationalized through the following discrete loss function, which penalizes deviation from exactness on the test basis:
\begin{equation}
   \mathcal C(\mathcal X;\theta) = \sum_{\boldsymbol j : \sum_{k=1}^d j_k \leq J}\left(\left(\sum_{i=1}^n
       W_i(\mathcal X ; \theta) \Phi_{\boldsymbol j}(x_i)\right) -1\right)^2.
\end{equation}
The training phase consists of multiple iterations where a new set of points $\mathcal X$ is obtained by resampling, and the model parameters $\theta$ are updated using the ADAM optimizer.

\section{Generating Non-Uniform Point Clouds}

To ensure robust training and evaluation of our neural network-based quadrature model, it is essential to be capable of generating point clouds that not only cover the entire domain but also feature non-uniform densities. This approach simulates the conditions commonly encountered in practical applications, such as adaptive mesh refinement or scenarios where physical phenomena are concentrated in specific areas. Our procedure for generating these point clouds is a two-step process that combines the superior uniformity of quasi-Monte Carlo sequences with a deterministic, non-linear transformation.

\subsection{Quasi-Random Base Sampling}

Instead of using entirely random numbers, which can lead to the
formation of gaps, we are generating sets of points with
low-discrepancy characteristics using a quasi-Monte Carlo (QMC)
method. QMC sequences are designed to cover the sampling space more
evenly than entirely random sampling strategies. Specifically, in this
work, we employ the Sobol' sequence\cite{sobol1967distribution}, a widely-used standard for QMC
integration.

Given $n$, we generate a base point set $\mathcal{S} = \{\mathbf{s}_1,
\mathbf{s}_2, \dots, \mathbf{s}_n\}$, where each point $\mathbf{s}_i$
is drawn from the Sobol' sequence within the unit hypercube $\Omega =
[0,1]^d$. This set $\mathcal{S}$ provides a basis for our sampling method.

\subsection{ Non-Linear Warping}

While the Sobol' sequence provides a well-behaved uniform distribution, our goal is to create point clouds with varying densities. To achieve this in a controllable and general manner, we can apply a bijective, non-linear transformation $\mathbf{G}: [0,1]^d \to [0,1]^d$ to uniformly distributed points. This vector-valued function warps the entire unit hypercube, allowing for complex, non-separable distortions beyond simple coordinate-wise stretching.

This transformation deterministically modifies the local density of points. Regions where the transformation compresses space lead to higher point density, while regions that are expanded result in lower density.

Each point $\mathbf{s}_i \in \mathcal{S}$ is to be mapped to a transformed point $\mathbf{x}_i$ to create the final non-uniform point cloud $\mathcal{X} = \{\mathbf{x}_1, \mathbf{x}_2, \dots, \mathbf{x}_n\}$. The transformation is defined as follows:

\begin{equation}
\mathbf{x}_i = \mathbf{G}(\mathbf{s}_i), \quad \text{for } i=1, \dots, n.
\end{equation}

This two-step procedure enables us to create on demand non-uniform point clouds with certain characteristics by simply modifying the warping function $\mathbf{G}(\cdot)$. The resulting unordered set $\mathcal{X}$ serves then as the direct input for our quadrature model, creating a robust framework to learn quadrature weights for arbitrary point clouds.

{
While domain transformations generate non-uniform point clouds, they can complicate numerical integration. Previous research indicates that warping quadrature rules obtained from a reference domain often degrade integration accuracy, particularly under complex mappings \cite{jakeman2019polynomial, warburton2013low, klebanov2023transporting, ernst2025learning}. On the other hand, \texttt{QuadrANN} learns to adapt to geometric distortions by adjusting weights based on local density features observed in the target point cloud.
}

\section{The QuadrANN Architecture}

To address the challenge of learning quadrature rules directly from
point clouds, we introduce a new neural network architecture,
\texttt{QuadrANN}. It is a deep Graph Neural Network (GNN)
architecture~\cite{scarselli2009graph, wu2021comprehensive} that
builds on principles from  graph models for point clouds
\cite{wang2019dynamic} and is designed to capture and utilize
geometric information at multiple scales effectively. Our approach is based on the principle that the
quadrature weight for a given point depends not only on its position
but also on its role within the geometry and structure of the
domain. The architecture is composed of three primary stages: an
initial geometric encoding layer, a number of globally-aware
propagation layers, and a final weight prediction neural network.

Several critical architectural decisions were made to optimize \texttt{QuadrANN} for this task. First, to capture rich local information, we combine two powerful features: Positional Encoding (PE), a technique popularized by Transformers and coordinate-based networks~\cite{vaswani2017attention, mildenhall2020nerf}, provides a high-frequency understanding of relative positions, while an explicit density estimate provides a direct, robust measure of local point concentration. Second, to ensure the entire shape informs the learned weights, each message-passing step is conditioned on a global context vector, a strategy used in seminal point cloud models~\cite{qi2017pointnet}. Third, inspired by \texttt{DenseNets}~\cite{huang2017densely}, we employ a dense concatenation of feature maps from all layers. This gives the prediction network simultaneous access to multi-scale features, mitigating the over-smoothing problem, common in deep GNNs~\cite{li2018deeper}.

Ultimately, the complete model network is trained using a global softmax function in the final layer. This activation function guarantees that the output weights are positive and sum to one, which satisfies the mathematical properties of a quadrature rule. The following sections will provide a detailed explanation of each component of this architecture.

The complete formulation of the model is presented in Algorithm~\ref{alg:quadrann}.

\subsection{Input Feature Engineering}

In order to prepare the raw point cloud, we perform an initial
feature engineering on the input coordinates. The raw coordinates
$\{\mathbf{x}_i\}_{i=1}^n$, assumed to lie within the unit hypercube
$[0, 1]^d$ (or in a subset of this domain). We first normalize the raw coordinates to the typical domain $[-1, 1]^d$ via the linear transformation $\mathbf{x}_i \leftarrow 2\mathbf{x}_i - 1$.

To enable the network to learn high-frequency patterns and to improve its ability to reason about both absolute and relative positions of the points, we employ Positional Encoding (PE). Using PE is a common technique that maps a low-dimensional continuous space into a higher-dimensional feature space by employing a set of sinusoidal functions with various frequencies.

For a general input vector $\mathbf{x} \in \mathbb{R}^d$, the positional encoding $\text{PE}(\mathbf{x})$ is constructed by applying the encoding to each component of the vector and concatenating the results:
\begin{equation}
\text{PE}_p(\mathbf{x}) = (\gamma(x_1;p), \gamma(x_2;p), \dots, \gamma(x_d;p))
\end{equation}
The encoding function $\gamma(\cdot;p)$ for a single scalar component $\alpha$ is defined as:
$$
\gamma(\alpha;p) = \left( \sin(2^0 \pi \alpha), \cos(2^0 \pi \alpha), \dots, \sin(2^{p-1} \pi \alpha), \cos(2^{p-1} \pi \alpha) \right)
$$
where $p$ is a hyperparameter representing the number of frequency
bands.
This multi-scale representation provides the model with a rich, unambiguous signal for both fine-grained detail and global position, resulting in a feature vector of dimension $d \times 2p$.

\subsection{Geometric Encoding Layer}

The first layer of our network creates a high-dimensional latent representation for each point $\mathbf{x}_i \in \mathbb{R}^d$. The dimensionality of these representations is denoted by $Q_1$, which we constrain to be an even number. An even number ensures that the feature dimension can be cleanly halved, as the final node features are formed by concatenating two vectors of equal size ($Q_1/2$). The resulting latent vector is denoted as $\mathbf{o}_i^{(1)} \in \mathbb{R}^{Q_1}$.

We begin by defining two $k$-NN graphs. First, a graph with $k'$ neighbours is used to compute a local density feature, $\rho_i \in \mathbb{R}$, for each point. The density feature is defined to be the reciprocal of the mean distance of a point's neighbours to their local centroid, providing a robust measure of point concentration. Specifically, for each point $\mathbf{x}_i$, we define its neighborhood centroid as $\mathbf{c}_i = \frac{1}{k'} \sum_{j:\mathbf{x}_j \in V_{k'}(\mathbf{x}_i)} \mathbf{x}_j$, and the density is then $\rho_i =
\left( \frac{1}{k'} \sum_{j:\mathbf{x}_j \in V_{k'}(\mathbf{x}_i)} \|\mathbf{x}_j -
  \mathbf{c}_i\|_2\right)^{-1}$. For computational efficiency, we
can always set $k' = k$, requiring the $k$-NN search to be performed
only once for each point $\mathbf{x}_i$.

By incorporating the explicit density feature, $\rho_j$, into the
message-passing scheme, the model acquires essential information about
local point concentration. Although a GNN can implicitly deduce
density from its inputs, supplying a precise measure enables the model
to reason about point density directly and frees up model capacity,
enabling the subsequent MLP to learn the more complex relationships
between geometry, density, and the resulting quadrature
weights. Ultimately, this leads to a more robust and accurate
quadrature rule.

Subsequently, a graph with $k$ neighbours determines the neighbourhood $V_k(\mathbf{x}_i)$, which facilitates the main message passing among the points. Within this neighbourhood, a message vector, $\mathcal{G}_{ij}$, is generated by augmenting the geometric features with the pre-computed density feature of the neighbours:
\begin{equation}
\mathcal{G}_{ij} = \left[ \mathbf{x}_j, \text{PE}_{p_1}(\mathbf{x}_j), \mathbf{x}_j - \mathbf{x}_i, \text{PE}_{p_2}(\mathbf{x}_j - \mathbf{x}_i), \rho_j \right]
\end{equation}
This message combines five key features:
\begin{enumerate}
\item the \textbf{absolute position} of the neighbor ($\mathbf{x}_j$) and its \textbf{positional encoding} ($\text{PE}_{p_1}(\mathbf{x}_j)$) to provide global context;
\item the \textbf{relative position vector} ($\mathbf{x}_j - \mathbf{x}_i$) to capture their Euclidean relationship;
\item the \textbf{positional encoding of the relative vector} ($\text{PE}_{p_2}(\mathbf{x}_j - \mathbf{x}_i)$) to provide a high-frequency representation of the local relationship;
\item the \textbf{local density explicit feature} ($\rho_j$) at the neighboring point to provide explicit information about the local concentration of the points.
\end{enumerate}

A dedicated MLP then processes this augmented data, denoted $\text{MLP}_1$. This network consists of three layers and its role is to transform the rich input data into a meaningful latent feature space:
\begin{enumerate}
    \item a linear transformation that maps the input $\mathcal{G}_{ij}$ to a hidden representation, followed by a GELU activation function;
    \item a second linear transformation and GELU activation further process the features;
    \item a final linear transformation that produces the output message $\mathbf{m}_{ij}^{(1)} \in \mathbb{R}^{Q_1/2}$.
\end{enumerate}
Finally, all learned messages are aggregated using a higher-order statistical function to form the initial node feature $\mathbf{o}_i^{(1)} \in \mathbb{R}^{Q_1}$:
\begin{equation}
\mathbf{o}_i^{(1)} = \left[ \boldsymbol{\mu}_i(\mathbf{m}_{ij}^{(1)}), \boldsymbol{\sigma}_i(\mathbf{m}_{ij}^{(1)}) \right]
\end{equation}
This initial encoding provides a rich, geometrically-aware basis for subsequent layers.

\subsection{Globally-Aware Propagation Layers}

Following the initial encoding, a stack of $L-1$ propagation layers refines the node features. At each layer $l \in \{2, \dots, L\}$, a global context vector $\mathbf{g}^{(l-1)} \in \mathbb{R}^{Q_1}$ is first computed by averaging all node features from the previous layer. A message is then computed by a layer-specific $\text{MLP}_l$, which has a compact two-layer structure sufficient for operating on already-learned features:
\begin{enumerate}
    \item a linear tranformation maps the concatenated input vector $[\mathbf{o}_j^{(l-1)}, \mathbf{o}_j^{(l-1)} - \mathbf{o}_i^{(l-1)}, \mathbf{g}^{(l-1)}]$ of dimension $3Q_1$ to a hidden representation, followed by a GELU activation function;
    \item a final linear transformation produces the output message $\mathbf{m}_{ij}^{(l)} \in \mathbb{R}^{Q_1/2}$.
\end{enumerate}
The new node features, $\mathbf{o}_i^{(l)} \in \mathbb{R}^{Q_1}$, are then formed by concatenating the mean and standard deviation of the messages. After the final propagation layer, we create a comprehensive representation, $\mathbf{o}_i^{\text{final}} \in \mathbb{R}^{L \cdot Q_1}$, by concatenating the feature vectors from all layers:
\begin{equation}
\mathbf{o}_i^{\text{final}} = \left[ \mathbf{o}_i^{(1)}, \mathbf{o}_i^{(2)}, \dots, \mathbf{o}_i^{(L)} \right].
\end{equation}

\subsection{Weight Prediction Network}

The final stage of the model maps the comprehensive feature vector
$\mathbf{o}_i^{\text{final}}$ to a quadrature weight $w_i$. This is
performed by a point-wise MLP, $\text{MLP}_{\text{out}}$, whose
structure is defined by a second hyperparameter, $Q_2$. We require
again $Q_2$ to be an even positive integer. This network is designed as follows:
\begin{enumerate}
    \item an input linear transformation takes the $L \times Q_1$-dimensional vector $\mathbf{o}_i^{\text{final}}$;
    \item a first hidden linear layer maps the input to a $Q_2$-dimensional
      space, followed by a GELU activation and a dropout (\(=0.5\)) layer for regularization;
    \item a second hidden layer reduces the dimension from $Q_2$ to
      $Q_2/2$, followed again by GELU and dropout (\(=0.5\));
    \item an output layer maps the final hidden representation to a single scalar logit, $v_i$.
\end{enumerate}

The vector of all logits, $\mathbf{v} = (v_1, \dots, v_n)^T$, is then normalized using a global \texttt{softmax} function to ensure the weights are positive and sum to unity ($w_i > 0, \sum_i w_i = 1$).

\begin{algorithm}[t]
\caption{The \texttt{QuadrANN} Architecture}
\label{alg:quadrann}
\begin{algorithmic}[1]
    \State \textbf{Input:} Point cloud $\mathcal{X} = \{\mathbf{x}_i\}_{i=1}^n \subseteq
    \mathbb{R}^{d}$, density neighbors $k'$, message neighbors $k$, layers $L$,
    PE bands $p_1, p_2$, node dimension $Q_1$, prediction network parameter $Q_2$.
    \Statex \textbf{Learnable Parameters:} Geometric MLP $\text{MLP}_1$, Propagation MLPs $\{\text{MLP}_l\}_{l=2}^L$, Prediction MLP $\text{MLP}_{\text{out}}$.
    \State \textbf{Output:} Quadrature weights $\mathbf{w} \in \mathbb{R}^{n}$.

    \Procedure{QuadrANN}{$\{\mathbf{x}_i\}_{i=1}^n$}
    \State \Comment{\textbf{Stage 1: Geometric Encoding Layer}}
    \State Find density neighborhood $V_{k'}(\mathbf{x}_i)$ for each point $\mathbf{x}_i$.
    \For{each point $\mathbf{x}_i,\ i = 1, \dots, n$} \Comment{Pre-compute density feature}
        \State $\mathbf{c}_i \gets \frac{1}{k'} \sum_{j:\mathbf{x}_j \in V_{k'}(\mathbf{x}_i)} \mathbf{x}_j$ \Comment{Neighborhood centroid}
        \State $\rho_i \gets ( \frac{1}{k'} \sum_{j:\mathbf{x}_j \in V_{k'}(\mathbf{x}_i)}
        \|\mathbf{x}_j - \mathbf{c}_i\|_2)^{-1}$ \Comment{Density
          feature is the inverse of mean distance}
    \EndFor
    \Statex
    \State Find message passing neighborhood $V_k(\mathbf{x}_i)$ for each point $\mathbf{x}_i$. \Comment{Can be skipped if $k'=k$}
    \For{each point $\mathbf{x}_i,\ i = 1, \dots, n$}
        \For{each neighbor $\mathbf{x}_j \in V_k(\mathbf{x}_i)$}
            \State $\mathcal{G}_{ij} \gets \left[ \mathbf{x}_j, \text{PE}_{p_1}(\mathbf{x}_j), \mathbf{x}_j - \mathbf{x}_i, \text{PE}_{p_2}(\mathbf{x}_j - \mathbf{x}_i), \rho_j \right]$ 
            \State $\mathbf{m}_{ij}^{(1)} \gets \text{MLP}_1(\mathcal{G}_{ij})$
        \EndFor
        \State $\boldsymbol{\mu}_i^{(1)} \gets \underset{j:\mathbf{x}_j \in V_k(\mathbf{x}_i)}{\text{Mean}}(\{\mathbf{m}_{ij}^{(1)}\})$; \quad $\boldsymbol{\sigma}_i^{(1)} \gets \underset{j:\mathbf{x}_j \in V_k(\mathbf{x}_i)}{\text{Std}}(\{\mathbf{m}_{ij}^{(1)}\})$
        \State $\mathbf{o}_i^{(1)} \gets [\boldsymbol{\mu}_i^{(1)}, \boldsymbol{\sigma}_i^{(1)}]$
    \EndFor
    \Statex

    \State \Comment{\textbf{Stage 2: Globally-Aware Propagation Layers}}
    \For{$l = 2, \dots, L$}
        \State $\mathbf{g}^{(l-1)} \gets \frac{1}{n} \sum_{p=1}^n \mathbf{o}_p^{(l-1)}$
        \For{each $i = 1, \dots, n$}
            \For{each $j:\mathbf{x}_j \in V_k(\mathbf{x}_i)$}
                \State $\mathbf{m}_{ij}^{(l)} \gets \text{MLP}_l([\mathbf{o}_j^{(l-1)}, \mathbf{o}_j^{(l-1)} - \mathbf{o}_i^{(l-1)}, \mathbf{g}^{(l-1)}])$
            \EndFor
            \State $\mathbf{o}_i^{(l)} \gets [\underset{j:\mathbf{x}_j \in V_k(\mathbf{x}_i)}{\text{Mean}}(\{\mathbf{m}_{ij}^{(l)}\}), \underset{j:\mathbf{x}_j \in V_k(\mathbf{x}_i)}{\text{Std}}(\{\mathbf{m}_{ij}^{(l)}\})]$
        \EndFor
    \EndFor
    \Statex

    \State \Comment{\textbf{Stage 3: Final Representation and Weight Prediction}}
    \For{each $i = 1, \dots, n$}
        \State $\mathbf{o}_i^{\text{final}} \gets [\mathbf{o}_i^{(1)}, \mathbf{o}_i^{(2)}, \dots, \mathbf{o}_i^{(L)}]$
        \State $v_i \gets \text{MLP}_{\text{out}}(\mathbf{o}_i^{\text{final}})$
    \EndFor
    \State $\mathbf{w} \gets \text{Softmax}((v_1, \dots, v_n)^T)$
    \State
    \State \Return $\mathbf{w}$
    \EndProcedure
\end{algorithmic}
\end{algorithm}

{
It is worth distinguishing the complementary roles of Positional Encoding (PE) and the Message Passing layers. While PE enables the network to resolve high-frequency spatial dependencies and distinguish unique point locations (addressing the spectral bias of standard MLPs), it remains a point-wise feature that does not explicitly encode point interactions. The message-passing layers are necessary to capture the local density of the point cloud, specifically how the warping function compresses or expands the space around a point.
This proposed architecture was determined through extensive testing, in which we evaluated the model by retaining only one component at a time. We observed that relying solely on PE or removing the local aggregation consistently degraded performance, particularly in high-dimensional scenarios. The proposed combination was thus identified as the necessary configuration to robustly compensate for the geometric distortions introduced by non-linear warping.

}

\section{Applications}

In this section, we demonstrate the robustness of the \texttt{QuadrANN} framework through a number of numerical experiments. We begin by directly evaluating its performance on pure integration problems and then showcase its utility as a component within complex variational solvers for time-dependent Partial Differential Equations (PDEs).

A challenging point in numerical integration, particularly for functions with sharp gradients or localised peaks, is the fact that standard sampling methods may not adequately cover the regions of interest. Methods like Monte Carlo or even quasi-Monte Carlo (QMC) that typically distribute points uniformly may undersample those critical parts of the domain where the contribution of the integrand is most significant, leading to slow convergence or high variance approximations. To achieve high accuracy with a limited number of points, the sampling distribution should ideally be concentrated in these important regions. Our approach is designed to excel in precisely this scenario by learning the correct quadrature weights for a given, potentially non-uniform, point distribution.

In order to ensure the robustness of our model, we test it in the context of non-uniform point clouds as described in the previous sections. These clouds are generated by first creating a low-discrepancy point set using a scrambled Sobol' sequence and then applying a non-linear, coordinate-wise warping function, $\mathbf{G}(\mathbf{x}) = (g(x_1), \dots, g(x_d))$. The function $g: [0,1] \to [0,1]$ is defined as:
\begin{equation}
 g(s) = 0.95s + 0.05\left(4(s-0.5)^3 + 0.5\right).
\end{equation}
This transformation intentionally concentrates points around the
center of the domain, simulating a scenario where the energy of the
integrand is localized around that point. By training
\texttt{QuadrANN} on a constantly changing set of warped point clouds,
the network learns a quadrature rule that is not tied to any single
cloud but is instead a general-purpose function adaptable to local
point density variations. Note that this exact transformation is used for all examples that follow.

Our numerical examples are divided into two parts. First, we present numerical integration examples to directly benchmark the accuracy and stability of \texttt{QuadrANN} against standard Sobol'-QMC. These tests include a simple 2D unit square, a non-convex 2D L-shaped domain, and a 4D hypercube to test the performance of the model under the curse of dimensionality. Second, we apply \texttt{QuadrANN} to solve PDEs: the Heat equation in 2D and the Fokker-Planck equation in 4D. These examples demonstrate the practical advantage of using an adaptive, data-driven quadrature rule within variational, time-stepping schemes that are common in modern scientific computing.

\subsection{Numerical Integration}

To test the performance of \texttt{QuadrANN}, we conducted
a number of numerical examples. For each example, we trained a
dedicated model to learn a quadrature rule for a specific domain and
according to a certain point cloud generation policy. In particular,
at each training epoch, we generate a batch of point clouds from a
scrambled Sobol' sequence and then apply the coordinate-wise warping
function $g(s)$ to create non-uniform distributions. The network is
updated based on these generated, non-uniform clouds. This ensures
that the model learns a general quadrature weight estimator adaptable
to point density variations, rather than overfitting to a single point
cloud. 

{ The performance of the model is evaluated over a large number of independent trials ($M=2^{12}$ in 2D, $M=2^{13}$ in 4D).}
In each trial, a new,
``unseen'' point cloud is generated, and \texttt{QuadrANN} is
used to compute the integral.

The final results are the
mean and standard deviation of these $M$ integral
estimates. For a straightforward comparison, the standard quasi-Monte Carlo
(QMC) baseline computes its  estimates on the corresponding $M$
original (unwarped) Sobol' point clouds considering uniform
weights.
In all cases, the objective is to integrate a known function, allowing
for explicit error analysis. The actual value of the integral is
$1$ for the first and third examples, and $3/4$ for the
second one.

\subsubsection{\underline{Example 1: 2D Integration on a Unit Square}}

The first example serves as a reference for performance in a simple domain, $\Omega = [0,1]^2$. The integrand here is a 2D Gaussian probability density function (PDF) selected for its smoothness and localized characteristics:
\begin{equation}\label{eq:integrand_f1}
f_1(\mathbf{x}) = \mathcal{N}(\mathbf{x} \,|\, \boldsymbol\mu, \sigma^2\mathbf{I}).
\end{equation}
For this example, the parameters are set to $\boldsymbol\mu=(0.5,
0.5)^T$ and $\sigma=0.025$ to center the function in the domain. The
\texttt{QuadrANN} model was configured with $L=3$ layers, feature
dimensions of $Q_1=2^7$, $Q_2=2^8$, $k=k^\prime=16$ nearest neighbors,positional encoding bands $p_1=3, p_2=5$, and $K_{\text{max}}=10$.

{
The results of this test case are summarized as follows. In this extensive trial
($M=2^{12}$), the \texttt{QuadrANN} model yielded a mean integral of $1.0033$ with a
standard deviation of $0.0874$. These results compare to the Sobol'--QMC method,
which resulted in a mean of $0.9991$ and a standard deviation of $0.1001$.
While Sobol' shows slightly lower bias in this instance, our model demonstrated
superior stability by attaining about $12.7$\% reduction in standard deviation over the
Sobol'--QMC.
}

\begin{figure}[h!]
    \centering
    \includegraphics[width=0.75\textwidth]{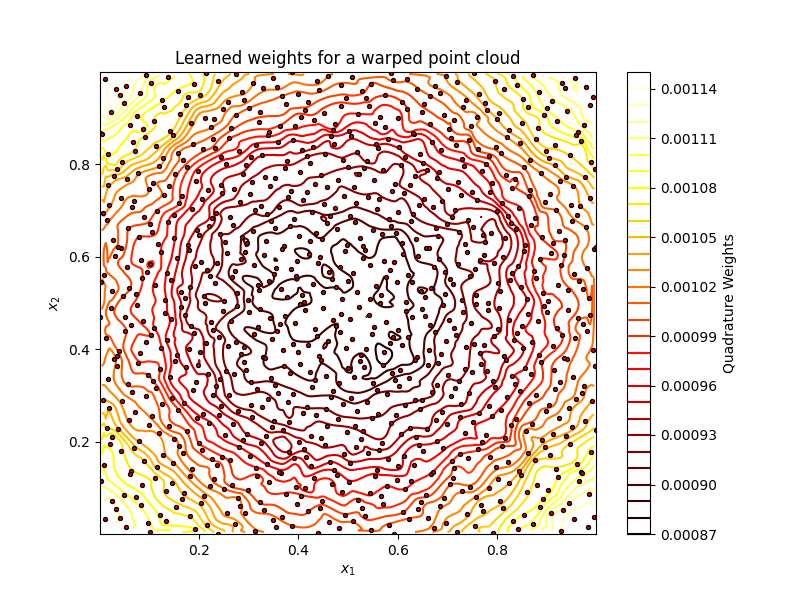}
    \caption{Interpolated quadrature weights learned by \texttt{QuadrANN} for the 2D unit square example. The color map indicates the weight value, and the dots refer to the locations of the points.}
    \label{fig:square_weights}
\end{figure}

Figure~\ref{fig:square_weights} visualizes the
learned quadrature weights corresponding to a random ``unseen'' point
cloud.

\subsubsection{\underline{Example 2: 2D Integration on a Non-Convex Domain}}\label{ex2:lshaped}

To challenge the model with complex geometry, the second example is conducted on a non-convex, L-shaped domain, defined as $\Omega = [0,1]^2 \setminus [0.5, 1]^2$. The inner corner of this domain makes it particularly difficult for standard methods. The integrand is the same Gaussian PDF defined in Equation~\eqref{eq:integrand_f1}.

This setup tests the capacity of \texttt{QuadrANN} to adapt
its weights to complex boundary. For this complex geometry,
we increased the network depth to $L=4$. The
additional propagation layer provides the model with a larger
receptive field and a greater capacity for feature extraction, which
is necessary to resolve the intricate non-convex boundary. The other
model hyperparameters were set as in the previous case $Q_1=2^7$,
$Q_2=2^8$, $k=k^\prime=16$, $p_1=3, p_2=5$, and $K_{\text{max}}=10$. Each point cloud consists of $n=2^{10}=1024$ points.
It should be noted that in this test case, the base functions were
re-normalized so that their integral over the L-shaped domain
evaluates to one. 
The model is trained for $256$ epochs.
{
For this test case with $M=2^{12}$ realizations, \texttt{QuadrANN} achieved a mean integral value
of $0.7504$ with a standard deviation of $0.0908$. The Sobol'--QMC
produced a mean of $0.7492$ with a standard deviation of $0.0972$. Here, the \texttt{QuadrANN} model achieved an absolute error of $0.0004$, whereas the corresponding Sobol'--QMC error was $0.0008$. Additionally, our model exhibited improved precision, achieving approximately a $6.6$\% reduction in standard deviation when compared to the Sobol'--QMC approach.
}

\begin{figure}[h!]
    \centering
    \includegraphics[width=0.8\textwidth]{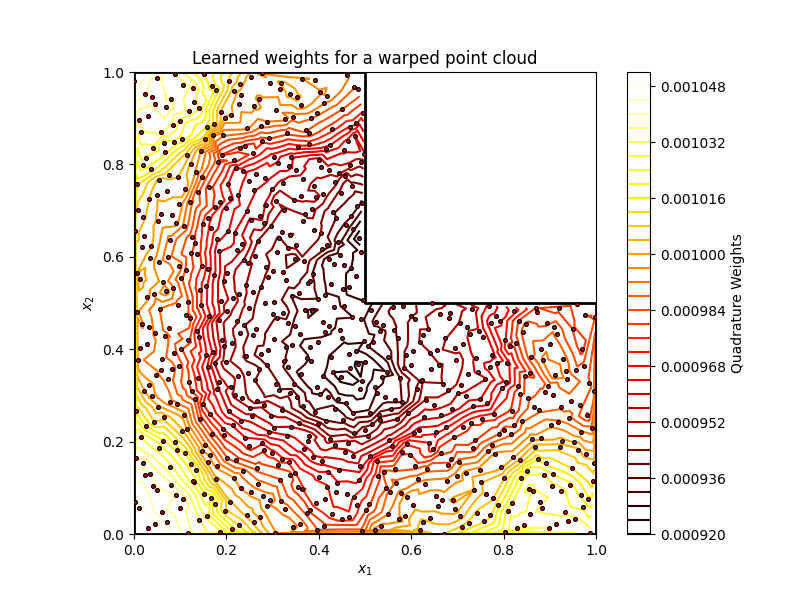}
    \caption{Visualization of the quadrature weights for the L-shaped
      domain. The model accurately assigns the weight of the points to
      adapt to the non-convex geometry. The color map indicates the
      weight value, and the dots refer to the locations of the points.}
    \label{fig:lshape_weights}
\end{figure}

The learned weights, shown in Figure~\ref{fig:lshape_weights},
demonstrate the ability of the model to adjust its weights for adapting to the non-convex geometry.

\subsubsection{\underline{Example 3: High-Dimensional Integration}}

The final pure integration example tests the performance of the model in a higher-dimensional space scenario, $\Omega = [0,1]^4$, where the curse of dimensionality is more pronounced. We compute the integral of a 4D Gaussian probability density function:
\begin{equation}
  f_3(\mathbf{x}) = \mathcal{N}(\mathbf{x} \,|\, \boldsymbol\mu, \sigma^2\mathbf{I}).
\end{equation}
with a mean $\boldsymbol\mu=(0.5,0.5,0.5,0.5)^T$, and a standard deviation $\sigma=0.075$. The network
for this 4D case used $L=4$ layers, feature dimensions of
$Q_1=2^7$, $Q_2=2^8$, $k=k^\prime=32$ neighbors, PE bands $p_1=5,
p_2=7$, and $K_{\text{max}}=6$. The model is trained for $512$ epochs,
% with a batch size of $4$,
using point clouds of size $n=2^{12}$.

{
While the preceding 2D examples primarily served to visualize weight
distribution, this 4D scenario highlights the sparsity challenges of
high-dimensional integration, necessitating a more rigorous
theoretical benchmark. To this end, we compare our method against
Kernel Optimal Quadrature (KOQ).
For our KOQ implementation, we employed a standard squared exponential (Gaussian) kernel, defined as $k(\mathbf{x}, \mathbf{y}) = \exp(-\|\mathbf{x} - \mathbf{y}\|_2^2 / 2s^2)$. To ensure the baseline remained robust across different random realizations of the point clouds without requiring manual tuning, the kernel bandwidth $s$ was determined adaptively using the median heuristic (i.e., setting $s$ to the median of the pairwise Euclidean distances of the point cloud). The optimal weights $\mathbf{w}$ were obtained by solving the linear system $(\mathbf{K} + \lambda \mathbf{I})\mathbf{w} = \mathbf{z}$, where $\lambda=10^{-5}$ is a regularization term, $\mathbf{K}$ is the Gram matrix of the kernel evaluations on the point cloud, and $\mathbf{z}$ is the vector of kernel mean embeddings with entries $z_i = \int_{[0,1]^4} k(x_i, \mathbf{y}) \, d\mathbf{y}$.

{
In this test (using $M=2^{13}$ realizations), \texttt{QuadrANN}
achieved a mean value of $1.0023$ and a standard deviation of
$0.0960$, outperforming standard Sobol'-QMC (mean $0.9986$, std
$0.1162$) and approaching the KOQ (mean $1.0009$, std
$0.0948$). However, a crucial distinction lies in computational
feasibility. While KOQ provides optimal weights for a static set, it
requires solving a dense linear system with $\mathcal{O}(n^3)$ complexity, which
is prohibitive for the dynamic resampling required in PDE solvers
(e.g., PINNs). \texttt{QuadrANN} bridges this gap, delivering
near-optimal accuracy comparable to KOQ but with a $\mathcal{O}(
\max(k,k^\prime )n^2)$ inference cost (dominated by the $k$-NN search), which is significantly faster than the cubic complexity of solver-based methods.
}

}

\subsubsection{\underline{Discussion on Stability}}

The results across all the examples underscore the benefits of our
approach, particularly in reducing the variance of the integral estimator. The consistent reduction in the standard deviation of the integral
estimates highlights a crucial feature of \texttt{QuadrANN}. While
these examples address pure integration problems, the implications for
extending to solve PDEs via functional minimization are clear.
The optimization process at the core of such solvers seeks the minimum of a discretized
integral. A quadrature rule with high variance provides a noisy
estimate of this functional, which can mislead the optimization
algorithm. By decreasing the variance, \texttt{QuadrANN} ensures a more
stable loss behaviour for the optimizer. This stability is critical
for ensuring the process converges robustly to the actual minimizer of
the functional and yields a more accurate final solution.

\subsection{Application to PDE Solvers}

We now demonstrate the practical efficacy of our ANN-based quadrature
rule within variational solvers for time-dependent PDEs. In each
scenario, the task is to find the solution by minimizing a
corresponding functional, where the accuracy of the numerical
integration is critical. Through these applications, we aim to
validate that a more accurate, data-driven quadrature rule leads to a
more faithful optimization of the underlying functional and,
consequently, a more precise solution.
The variational time-stepping approach used here is inspired by
minimizing movement schemes (IMEXs), which have proven effective in other complex deep learning contexts \cite{georgoulis2025deep}.

\subsubsection{\underline{PDE Solution Network and Setup}}

For the forthcoming PDE examples, the time-dependent solution $u(t_k,
\mathbf{x})$ is approximated by a neural network, $U(t_k, \mathbf{x};
\theta)$. We employ a DGM-type architecture \cite{sirignano2018dgm}, detailed in Algorithm~\ref{alg:dgm_forward}, which is proper for high-dimensional problems.

\begin{algorithm}[ht]
\caption{The DGM-type Solution Network Architecture}
\label{alg:dgm_forward}
\begin{algorithmic}
    \State \textbf{Input:} Spatial coordinates $\mathbf{x} \in \mathbb{R}^d$, time $t$.
    \Statex \textbf{Learnable Parameters:} $\{\mathbf{W}^{(\text{in})}, \mathbf{b}^{(\text{in})}, \{\mathbf{V}^{(s,l)}, \mathbf{W}^{(s,l)}, \mathbf{b}^{(s,l)}\}_{l=1}^L, \mathbf{W}^{(\text{out})}\}$ for $s \in \{g,z,r,h\}$.
    \Statex \textbf{Hyperparameters:} Number of layers $L$.
    \State \textbf{Output:} Approximate solution $U(t, \mathbf{x}; \theta)$.

    \Procedure{U}{$t, \mathbf{x}$}
        \State \Comment{Input Layer}
        \State $S^{(0)} \gets \tanh(\mathbf{W}^{(\text{in})} [t, \mathbf{x}]^T + \mathbf{b}^{(\text{in})})$ \Comment{Concatenate time and space}
        \State
        \State \Comment{DGM Hidden Layers}
        \For{$l = 1, \dots, L$}
            \State $G^{(l)} \gets \tanh(\mathbf{V}^{(g,l)} \mathbf{x} + \mathbf{W}^{(g,l)} S^{(l-1)} + \mathbf{b}^{(g,l)})$
            \State $Z^{(l)} \gets \tanh(\mathbf{V}^{(z,l)} \mathbf{x} + \mathbf{W}^{(z,l)} S^{(l-1)} + \mathbf{b}^{(z,l)})$
            \State $R^{(l)} \gets \tanh(\mathbf{V}^{(r,l)} \mathbf{x} + \mathbf{W}^{(r,l)} S^{(l-1)} + \mathbf{b}^{(r,l)})$
            \State $H^{(l)} \gets \tanh(\mathbf{V}^{(h,l)} \mathbf{x} + \mathbf{W}^{(h,l)}( S^{(l-1)} \odot R^{(l)}) + \mathbf{b}^{(h,l)})$ 
            \State $S^{(l)} \gets (1 - G^{(l)}) \odot H^{(l)} + Z^{(l)} \odot S^{(l-1)}$ \Comment{The new hidden state}
        \EndFor

        \State
        \State \Comment{Estimated solution}
        \State $U \gets \text{Softplus}(\mathbf{W}^{(\text{out})} S^{(L)})$
        \State
        \State \Return $U$
    \EndProcedure
\end{algorithmic}
\end{algorithm}

The solution network $U(t, \mathbf{x}; \theta)$ is configured with $L=3$
hidden layers and a width of $32$ nodes per layer. For the 2D Heat
Equation, we use a time step of $\tau = 0.01$. The optimization for the
first time step ($k=1$) is performed for $2^{11}$ epochs, while all
subsequent steps are trained for $2^7$ epochs with a batch size of
$8$. For the 4D Fokker-Planck equation, the time step is $\tau =
0.01$. The first step is trained for $2^{11}$ epochs and the rest for
$2^{7}$ epochs, using now larger point clouds of $n=2^{12}$ points, with
a batch size of $4$.

\subsubsection{\underline{Example: Heat Equation}}\label{example2:heat}
\begin{figure}[t]
    \centering
    \includegraphics[width=0.7\textwidth]{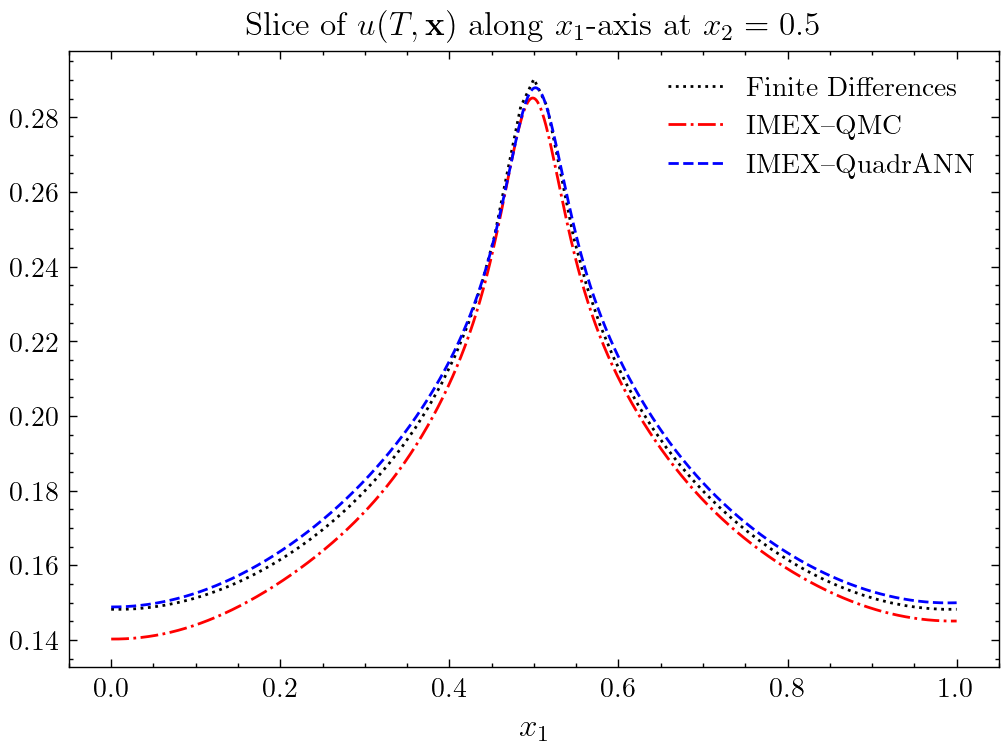}
    \caption{Comparison of solutions for the 2D Heat Equation. A 1D slice of the computed solution $u(T,\mathbf{x})$ at the final time $T=1$, taken along the $x_1$-axis at $x_2=0.5$. The plot compares the solution obtained using our proposed \texttt{QuadrANN} integration with the standard Sobol'-based Monte Carlo (MC) method against the Finite Difference Method (FDM) benchmark. The \texttt{QuadrANN}-based solution shows a closer agreement with the benchmark.}
    \label{fig:slice-heat}
  \end{figure}

\begin{figure}[t]
    \centering
    \includegraphics[width=1.\textwidth]{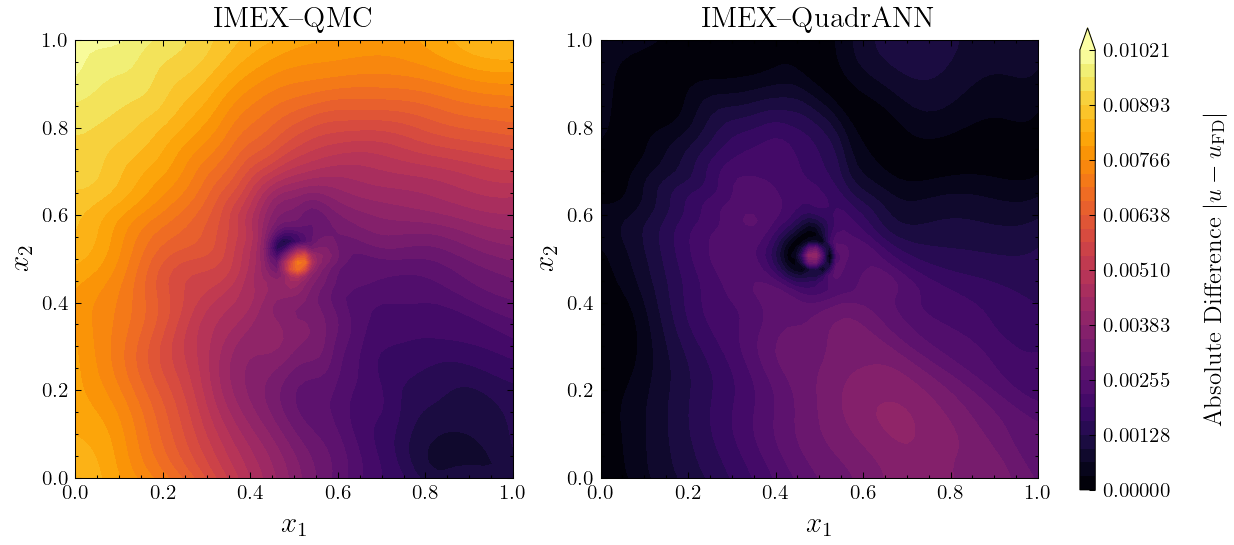}
\caption{Absolute error at final time \(T=1\) relative to the
      Finite Difference benchmark. The plots compare the error of the
      solution obtained with \texttt{QuadrANN} integration versus the
      corresponding Sobol'-QMC solution.}
    \label{fig:heat-colormap}
  \end{figure}

We wish to solve the heat equation on $\Omega = [0,1]^2$. The problem is to find $u(t,\mathbf{x})$ for $t\in [0,1]$ and $\mathbf{x} = (x_1,x_2) \in \Omega$, such that:
\begin{equation}
\begin{cases}
\partial_t u  - c \Delta u = f(t,\mathbf{x}) & \text{in } (0,1]\times \Omega,\ c = 1/8 \\
u(0,\mathbf{x}) = 0 & \text{in } \Omega \\
\frac{\partial u}{\partial \mathbf{n}}(t,\mathbf{x}) = 0 & \text{on } (0,1]\times \partial\Omega
\end{cases}
\end{equation}
The source term $f(t, \mathbf{x})$ is defined as:
\begin{equation}
f(t, \mathbf{x}) = 
\begin{cases} 
    50 \left(1 + \cos\left(\frac{\pi r(\mathbf{x})}{R}\right)\right) \left(2 + \cos(4\pi t)\right)e^{-3t} & \text{if } r(\mathbf{x}) < R \\
    0 & \text{if } r(\mathbf{x}) \ge R
\end{cases},
\end{equation}
where $r(\mathbf{x}) = \sqrt{(x_1-0.5)^2 + (x_2 - 0.5)^2}$ and $R=0.05$. We estimate the solution $u^k(\mathbf{x})$ at time $t_k=k\tau$ by minimizing the cost functional $\mathcal{J} [u^k]$:
\begin{equation}
\mathcal{J} [u^k]= \int_\Omega\frac{1}{2}\left(u^k(x) - u^{k-1}(x)\right)^2 + \frac{1}{2}\tau c \left(\nabla u^k(x)\right)^2 - \tau f(t^k,x)u^{k-1}(x) dx.
\end{equation}
The quadrature model is first trained on the $[0,1]^2$ domain for $256$ epochs,
% using a batch size of $32$,
where each batch consists of a point cloud of $2^{10}$ points created
using the adopted warping policy. The trained model is then used to
provide the quadrature weights for minimizing $\mathcal{J} [u^k]$. The solution
obtained using \texttt{QuadrANN} integration achieved a mean absolute
error of $1.527\times10^{-3}$, a clear improvement over the error of
$4.329\times10^{-3}$ from the solver using standard Sobol'-based Monte Carlo
integration. This confirms the advantage of our approach for about
15\% computational overhead. A qualitative comparison is presented in
Figure~\ref{fig:slice-heat}, { while the entire spatial distribution of the absolute error is depicted in Figure~\ref{fig:heat-colormap}.}

\subsubsection{\underline{Example: Fokker-Planck Equation}}\label{example3:fp}
To test our method in a higher-dimensional scenario, we focus on the
Fokker-Planck equation within the $\Omega = [0,1]^4$.

The problem describes the probability density function,
\(u(t,\textbf{x})\), for the state of particle whose position \(X_t\)
evolves according to an Ornstein-Uhlenbeck type Stochastic Differential Equation (SDE):
\begin{equation}
    d\mathbf{X}_t = -\gamma (\mathbf{X}_t - \mathbf{x}_0) dt + \sigma d\mathbf{W}_t,
\end{equation}
where the drift is $\gamma=0.5$, the mean-reversion level is $\mathbf{x}_0
= (0.5, 0.5, 0.5, 0.5)^T$, the diffusion is $\sigma=0.2$, and
$\mathbf{W}_t$ is a 4D Brownian motion.
We assume again Neumann boundary conditions on $\partial\Omega$.
The corresponding Fokker-Planck equation, which governs the time evolution of the probability density function $u(t,\mathbf{x})$, is given by:
\begin{equation}
  \partial_t u = \nabla \cdot [\gamma(\mathbf{x} - \mathbf{x}_0)u] + \frac{\sigma^2}{2} \Delta u.
\end{equation}
The considered initial condition, $u(0,\mathbf{x}) = u_0(\mathbf{x})$, is a mixture of two multivariate Gaussian distributions creating a bimodal landscape:
\begin{equation}
    u_0(\mathbf{x}) = 0.5 \ \mathcal{N}(\mathbf{x} | \boldsymbol\mu_1, \sigma_1^2 \mathbf I) + 0.5 \ \mathcal{N}(\mathbf{x} | \boldsymbol\mu_2, \sigma_2^2 \mathbf I),
\end{equation}
with means $\boldsymbol\mu_1 = (0.35, 0.35, 0.35, 0.35)^T$,
$\boldsymbol\mu_2 = (0.65, 0.65, 0.65, 0.65)^T$, and standard deviations of $0.1$ and $0.11$, respectively.

Following a variational time-stepping scheme, we approximate the solution $u^k(\mathbf{x})$ at each time step by minimizing the cost functional $\mathcal{J}[u^k]$:
\begin{equation}
\begin{split}
\mathcal{J}[u^k] = \int_\Omega \bigg[ & \frac{1}{2}\left(u^k(\mathbf{x}) - u^{k-1}(\mathbf{x})\right)^2 + \frac{\tau \sigma^2}{4} |\nabla u^k(\mathbf{x})|^2 \\
& - \tau \gamma \left((\mathbf{x}-\mathbf{x}_0) \cdot \nabla u^{k-1}(\mathbf{x})\right)u^k(\mathbf{x}) - 2\tau\gamma u^k(\mathbf{x}) \bigg] d\mathbf{x}.
\end{split}
\end{equation}
A benchmark solution is obtained by simulating the SDE for $10^7$
particles and using Kernel Density Estimation (KDE) with a bandwidth
parameter of $0.03$, a standard non-parametric technique \cite{hastie2009elements}, on the final particle positions to construct a high-fidelity density.
The error of our solver is quantified by the Hellinger distance to this benchmark.

\begin{figure}[t]
    \centering
    \includegraphics[width=0.7\textwidth]{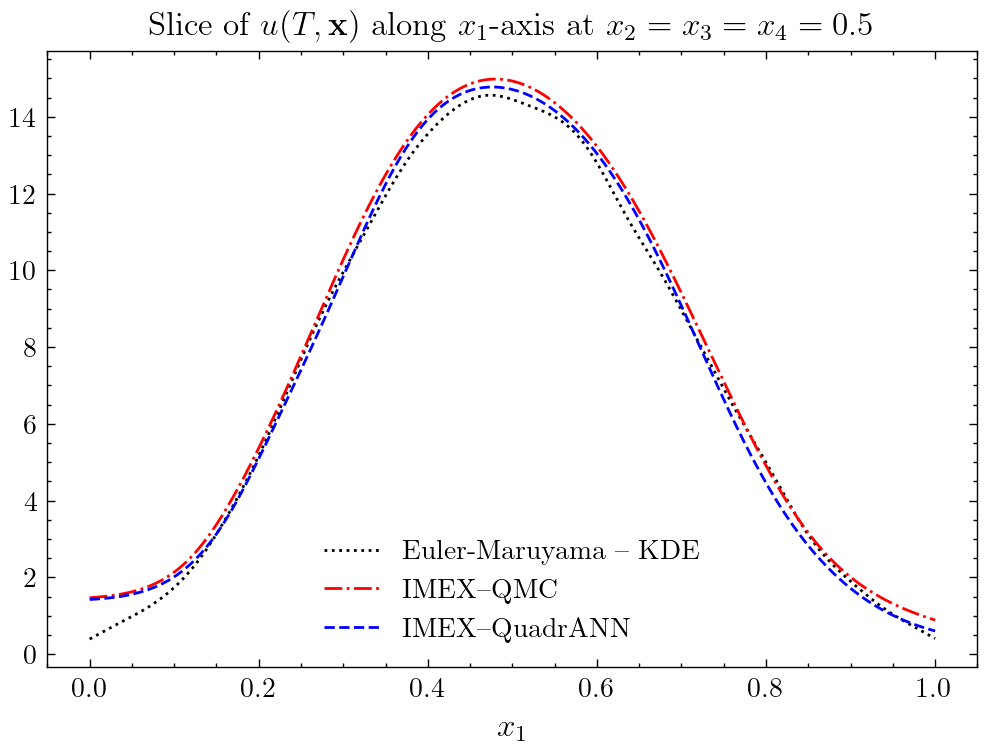}
    \caption{Comparison of solutions for the 4D Fokker-Planck equation. A 1D slice of the computed density $u(T,\mathbf{x})$ at the final time $T=1/2$, taken along the $x_1$-axis with other coordinates fixed at the center of the domain ($x_2=x_3=x_4=0.5$). The plot compares the density profile from our \texttt{QuadrANN} integration scheme with a standard Sobol'-based Monte Carlo (quasi-MC) method against the benchmark density constructed via SDE path simulation and Kernel Density Estimation (KDE).}
    \label{fig:slice-FP}
\end{figure}

For this 4D problem, a non-uniform point cloud is generated by
applying a coordinate-wise warping function $g(s)$ to a 4D Sobol'
sequence. A dedicated \texttt{QuadrANN} model is trained for $512$
epochs on point clouds of $n=2^{12}$. This
trained model then supplies the weights for the variational
solver. The final solution using \texttt{QuadrANN} integration
achieved a Hellinger distance of $0.0769$ to the KDE benchmark. This
is considered an improvement over the Hellinger distance of $0.0917$
obtained with standard Sobol'-QMC integration, achieved at a modest
20\% increase in computational time. Figure~\ref{fig:slice-FP}
provides a visual comparison of the obtained solutions.

\section{Conclusions}

In this work, we introduced \texttt{QuadrANN}, a deep learning architecture designed to overcome the critical challenge of numerical integration on non-uniform point clouds, a common task in modern mesh-free variational solvers for PDEs. Mainstream methods like quasi-Monte Carlo struggle to handle such distributions, while the traditional permutation-variant neural networks make them unsuitable for unordered point sets. Our approach treats these limitations by adopting a Graph Neural Network (GNN) to learn data-driven, permutation-invariant quadrature rules directly from geometry. The \texttt{QuadrANN} architecture introduces a rich geometric message-passing scheme that combines absolute and relative positional features as well as a density feature in its initial layer. In contrast, the subsequent layers exploit a global context vector to make the learned rule aware of both local point density and the overall domain shape.

We demonstrated the efficiency of our method through a series of numerical experiments. Direct integration tests on simple, non-convex, and high-dimensional domains confirmed that \texttt{QuadrANN} consistently and significantly reduces the variance of the integral estimate compared to standard Sobol'-QMC integration. The uncertainty reduction led to improved performance in solving PDE problems.
In particular, when the quadrature model was integrated into variational ML-solvers for the Heat and Fokker-Planck equations, our method yielded more accurate final solutions than those obtained using standard QMC quadrature, validating its practical merit.

The core benefit of this work is that it provides a stable framework for the optimization process used in mesh-free variational solvers. 
Particularly, \texttt{QuadrANN} ensures a more reliable convergence to
the energy minimizer by reducing the variance of the estimations of
the integral, even while providing mean estimates of the integrals comparable to those from Sobol'-QMC.

{
A rational extension of this work would be to integrate the
point-selection mechanism into the learning framework. Since picking
$n$ independent points yields a highly nonconvex landscape prone to
local minima, a more rational approach is to learn the warping
mechanism instead. Future research could explore parameterizing a
smooth, invertible mapping that transforms a standard low-discrepancy
sequence into an optimal point distribution tailored to the
integrand's features. By learning the transformation mechanism rather
than individual coordinates, we can achieve a fully adaptive
quadrature method that optimizes both the point cloud and its
associated weights, further enhancing efficiency for complex
variational minimization problems.
}

\section*{QuadrANN Implementation}
We have developed \texttt{QuadrANN} using \texttt{Pytorch-Geometric}\cite{pyg} on
top of \texttt{PyTorch}\cite{pytorch}.
Interested readers can refer to the following \texttt{GitHub}
repository (\href{https://github.com/kesmarag/QuadrANN}{github.com/kesmarag/QuadrANN}) for the source code.

\section*{Declarations}

\textbf{Funding:} Funding not applicable.

\printbibliography

\end{document}